\documentclass[12pt]{article}
\usepackage{amsmath,amssymb,bm}
\usepackage{color,url}
\usepackage{xcolor} 
\usepackage[numbers]{natbib} 
\usepackage{diagbox} 
\usepackage{multirow} 
\usepackage{dsfont} 
\newtheorem{algorithm}{Algorithm}

\def\re{{\mathbb{R}}}

\def\norm#1{\left\Vert{#1}\right\Vert}

\def\bee{\begin{equation}}
\def\ene{\end{equation}}
\def\beea{\begin{eqnarray}}
\def\enea{\end{eqnarray}}
\def\beeas{\begin{eqnarray*}}
\def\eneas{\end{eqnarray*}}

\def\mr#1{(\ref{#1})}

\def\RE{\hbox{Re}}

\def\ignore#1{}
\def\noise{{\bm{\xi}}}
\def\anoise{{\breve{\bm{\xi}}}}
\def\tnoise{{\breve{\bar{\bm\xi}}}}

\def\T{\text{T}}

\title { {\bf
Optimized convergence of stochastic gradient descent by weighted averaging }}

\author{\small Melinda Hagedorn, Heinrich Heine Univ., D\"usseldorf, Germany\\
  \small Florian Jarre, Heinrich Heine Univ., D\"usseldorf, Germany }
\date  { Sept. 23, 2022 \\ \phantom{nix} \ \ 
 \phantom{nix} \\
In memory of Oleg Burdakov}

\begin{document}
\maketitle 

\begin{abstract}
  Under mild assumptions 
  stochastic gradient methods asymptotically  achieve an optimal rate of convergence if
  the arithmetic mean of all iterates is returned as an approximate optimal solution.
  However, in the absence of stochastic noise, the arithmetic mean of all iterates converges
  considerably slower to the optimal solution than the iterates themselves.
  And also in the presence of noise, when a finite termination of the stochastic gradient method
  is considered, the arithmetic mean is not necessarily the best possible approximation
  to the unknown optimal solution.
  This paper aims at identifying optimal strategies in a particularly simple case,
  the minimization of a strongly convex function with i.\,i.\,d. noise terms and finite
  termination. Explicit formulas for the stochastic error and the optimization error are
  derived in dependence of certain parameters of the SGD method.
  The aim was to choose parameters such that both stochastic error and optimization error are
  reduced compared to arithmetic averaging. This aim could not be achieved; however, by allowing
  a slight increase of the stochastic error it was possible to select the parameters such
  that a significant reduction of the optimization error could be achieved. This reduction of
  the optimization error has a strong effect on the approximate solution generated by the
  stochastic gradient method in case that 
  only a moderate number of iterations is used or when the initial error is large. 
  The numerical examples confirm the theoretical results and suggest that a generalization
  to non-quadratic objective functions may be possible.
  
\ignore{
The aim of this paper is to unravel why averaging over all SGD iterates with the same weight leads asymptotically to the optimal rate of convergence. The approach is to consider the SGD method as a perturbed gradient descent method with i.\,i.\,d. noise terms. Although it is not possible to reduce the optimization error and the stochastic error simultaneously, in case of a limited number of descent steps, moderately increasing weights can decrease the optimization error without significantly deteriorating the stochastic error. Nevertheless, if the number of descent steps tends to infinity, the stochastic effects dominate and thus determine the rate of convergence more than the condition number of the objective function. This explains why the choice of equal weights is asymptotically optimal. In the theoretical analysis the objective function is assumed to be smooth and strongly convex. The numerical examples confirm the theoretical results and suggest that a generalization to non-quadratic objective functions may be possible.}
\end{abstract}


\newpage
\section{Introduction}
In Polyak and Juditsky \cite{polyak} it is shown that the
stochastic gradient descent algorithm asymptotically
achieves optimal complexity
  if short but constant step lengths are chosen and if the
  average over all iterates is used as final output.
  The aim of the current paper is to explain this phenomenon with
  a slightly different approach and to optimize the
  results further by the use of weighted averages
  while considering finite termination and the nature of
  the function to be minimized.\\
The main source of this paper, Polyak and Juditsky \cite{polyak}, also inspired many other scientists to investigate weighted averages in conjunction with the stochastic gradient method in more detail. In Neu and Rosasco \cite{neu}, a variant of the weighted average SGD with geometrically decreasing weights is analyzed in context of linear least-squares regression.
Likewise, Cohen and Nedi$\acute{\text{c}}$ \cite{cohen} deal with the least-square regression. They consider constrained problems and derive upper bounds for the convergence rate and for the asymptotic ratio between convergence rate and empirical risk minimizer depending on the dimension. More abstractly, papers like Izmailov et al. \cite{izmailov} and Guo et al. \cite{guo}, which deal with stochastic weight averaging, also build on \cite{polyak}.\\
Recently Sebbouh et al. \cite{sebbouh} have shown almost sure convergence rates for weighted average SGD with decreasing step size. This result has already been supplemented for strongly-convex and non-convex objective functions by Liu and Yuan \cite{liu}.\\
There certainly is also research in this area independent of Polyak and Juditsky \cite{polyak}. For example, the sampling of Needell et al. \cite{needell} and the paper of Shamir and Zhang \cite{shamir}, in which the polynomial-decay averaging and the suffix averaging are examined.\\
There are numerous further modifications of the stochastic gradient approach such using
  momentum or heavy-ball iterations, see, e.g., \cite{loizou}, variance reduction \cite{johnson},
  stochastic gradient boosting \cite{friedman},
  and modifications tailored to specific applications.
  The above is merely a short and incomplete selection of related work. It seems, however, that
  the focus of the present paper has not been considered in this form before.
  This paper returns to a simple general format as considered in \cite{polyak},
  aims at optimizing two parameters in this approach, and considers the case
  of an infinitely large training set in the numerical examples.
  A brief outline is given next.

\subsection{Outline}
A consequence of the results by
Polyak and Juditsky \cite{polyak} is that asymptotically
the optimal rate of convergence of a stochastic gradient descent method
is obtained when averaging {\em all} iterates with the
same weight. This may seem counter-intuitive as one would expect the later
iterates to be closer to the optimal solution and would therefore
allow higher weights for later iterates.
A short and intuitive explanation of why
  averaging over all iterates is optimal is given by the observation that
  the square root of the function value is reduced at a linear rate when
  sufficiently short steps with constant step lengths are chosen
  for minimizing a smooth, strongly convex function, while the
  reduction of the variance has a much slower rate of convergence.
  Thus, asymptotically, the stochastic effects determine the overall
  rate of convergence, and not the condition number of the function
  to be minimized. And from a stochastic point of view, averaging
  over all iterates with the same weight is a simple but optimal strategy.
  However, if only a limited number of stochastic descent steps are
  taken, the optimization aspect and the aspect of stochastic convergence
  need to be balanced to each other.
  The present paper is an attempt to derive a simple strategy
  that does so in an optimized form.

To this end, an elementary
derivation of the optimality result in \cite{polyak}
is attempted in Section \ref{sec.opt} for a particularly simple
situation, the minimization of a strongly convex quadratic function
$f:\re^n\to\re$,
\begin{align}\label{f}
f(x)\ \equiv\ \frac 1m \sum_{i=1}^m f_i(x) \ \equiv \ \frac 1m
\sum_{i=1}^m \frac 12 x^\T A^{(i)} x+(b^{(i)})^\T x+c_i
\end{align}
by a stochastic gradient approach.
 Exact and computable formulas for the variance of the averaged
iterates as a measure of the {\em stochastic error}
and for the contraction constant of the descent steps
as a measure of the {\em optimization error} of the weighted iterates
are derived. Given these explicit formulas standard nonlinear
minimization algorithms with different starting points were
applied to reduce both errors simultaneously by adjusting
certain parameters associated with the weights and the step lengths
of a stochastic gradient method. While it turned out that the goal
could not be reached of reducing both, optimization error and stochastic
error at the same time, a significant reduction of the optimization error was
possible by allowing a slight increase (e.g.\ of 10\%) of the
stochastic error.

Generalizations to further classes of smooth convex functions are discussed
   in Section \ref{sec.num}.

\subsection{Notation}
The condition number with respect to the
2-norm of a square matrix $A$ is denoted by $\hbox{cond}(A)$ and
the smallest and largest eigenvalues are denoted
by $\lambda_{min}(A)$ and  $\lambda_{max}(A)$.
The Hadamard product (componentwise product) of two matrices $X,Y$
with the same dimensions is denoted by $X\circ Y$
and powers of vectors are also defined componentwise,
e.g.\ $x^2=x\circ x$.
The diagonal of a square matrix $A$ is denoted by $\hbox{diag}(A)$
and a diagonal matrix with diagonal $x\in\re^n$ is denoted by $\hbox{Diag}(x)$.
The $i$-th canonical unit vector is denoted by $e_i$
and the all-one-vector is denoted by $e$, its dimension being given
by the context. Finally, let $\mathds{1}_{n\times n}$ denote the $n\times n$ matrix
with all entries equal to one.

\subsection{A Standard Stochastic Descent Method}
For large values of $m$ (or when $m=\infty$)
  a stochastic gradient descent method of the following form is considered:
Assume that a batch $S_k$  is chosen i.\,i.\,d. from the uniform
distribution of $\{1,\ldots,m\}$. Then, as is
well known, the expected value of the gradient of
$f_{S_k} (x^k):= \frac 1{|S_k|} \sum_{i\in S_k} f_i(x^k)$
is the full gradient, 
$$
 E(\nabla f_{S_k} (x^k)) = E\left(\frac 1{|S_k|}\sum_{i\in S_k}\nabla f_i(x^k)\right) 
= \nabla f(x^k).
$$
This motivates the stochastic gradient descent that uses the approximation $\nabla f_{S_k} (x^k)$ instead of $\nabla f(x^k)$
to define a sequence of iterates.
Let $\gamma: \re_+\to \re_+$ be a weakly monotonously decreasing function
  (used in the analysis below)
  and consider a-priori-defined\footnote{%
    In the presence of noise, a line search is difficult.}
  step lengths $\gamma_k=\gamma(k)$ in a stochastic gradient
descent approach
\bee \label{proc1}
x^{k+1} = x^k-\gamma_k\nabla f_{S_k} (x^k)= x^{k}- \gamma_{k} (\nabla f(x^{k})+\noise^{k})
\ene
with random noise term
\begin{align}\label{noisek}
\noise^k\  :=\ \nabla  f_{S_k} (x^k)-\nabla f(x^k) 
\ = \ \frac 1{|S_k|} \sum_{i\in S_k}  A^{(i)} x^k+b^{(i)}-
\frac 1m \sum_{i=1}^m  A^{(i)} x^k+b^{(i)}
\end{align}
whose expected value is zero.

\subsubsection{Assumptions}
For simplicity
it is assumed that the batch size $|S_k|$ is constant for all $k$.
As $\noise^k$ depends on the current iterate, and thus on previous noise terms
$\noise^i$ for $i<k$
it may seem unrealistic to assume that the terms $\noise^k$ are i.i.d.\ as well.
However, 
setting
$$
\Delta A_{S_k} :=\frac 1{|S_k|} \sum_{i\in S_k}  A^{(i)} -
\frac 1m \sum_{i=1}^m  A^{(i)} \ \ \hbox{ and } \ \ 
\Delta b_{S_k} :=\frac 1{|S_k|} \sum_{i\in S_k}b^{(i)}-\frac 1m
\sum_{i=1}^m  b^{(i)}
$$
the noise can be written as
$$
\noise^k =\Delta A_{S_k}  x^k+\Delta b_{S_k} 
$$
where both $\Delta A_{S_k} $ and $\Delta b_{S_k} $
are i.i.d. and the expected values satisfy
$$
E(\Delta A_{S_k} ) = 0\in\re^{n\times n} \qquad \hbox{and}  \qquad
E(\Delta b_{S_k} )=0\in\re^n.
$$
Since $|S_k|$ is constant for all $k$, also
$\sigma_A^2:=E(\|\Delta A_{S_k} \|_F^2)$ and
$\sigma_b^2:=E(\|\Delta b_{S_k} \|_2^2)$
are independent of $k$.
Moreover, when the iterates $x^k$ converge to some limit $x^*$
then, asymptotically, the
error terms
\bee\label{iid}
\noise^k \quad \hbox{are i.i.d.}
\ene
This is the situation that is also analyzed in a different context as part (c)
of Theorem 1 in \cite{polyak}, and this assumption will be used
for simplicity below.

For the analysis it can be assumed further, without loss of generality,
that $x^*=0$
and that the Hessian of $f$ is a diagonal matrix $D$ with diagonal elements
$0<D_{1,1} \le \ldots \le D_{n,n}$. Thus, $\nabla f(x)=Dx$.

Finally, assume that
some upper bound for the largest eigenvalue $D_{n,n}$ is known so that 
the step lengths $\gamma_{k}$ can be chosen in the half-open
interval $(0,1/ D_{n,n}]$. (A positive lower bound for $D_{1,1}$
  is not assumed to be known.)

Up to a factor of 2, the bound on $\gamma_k$
is essentially the bound in \cite{polyak}.
(Note that there is a minor error in \cite{polyak}:
in Assumption 2.2 (page 839) they write $2/\min_i \RE(D_{i,i})$ while
the correct statement would be
$\min_i 2 \RE(D_{i,i})/|D_{i,i}|^2$
in order for the
argument at the bottom of page 844 of \cite{polyak} to be valid.)

\subsection{A Simple Algorithm with Two Parameters}
In Polayk and Juditsky \cite{polyak}, an algorithm
  is considered with short step lengths and with output given by the arithmetic
  mean of all iterates $x^k$ generated via \mr{proc1}. 
  In the absence of noise,
  i.e.\ in case that all terms $\noise^k$ are zero, it is clear
  that the final iterate $x^k$ is a much better approximation
  to the optimal solution than the average.
  On the other hand, as shown in \cite{polyak},
  averaging all iterates with equal weights
  reduces the variance, and is asymptotically optimal for large $k$.
  In the following a method is considered that uses weighted averages
  with higher weights on later iterates. When the noise is small,
  the weighted average also is a better approximation to the optimal
  solution than the average with equal weights.
  To reduce the variance of the iterates with higher weight,
  a possible step length reduction for the later iterates is
  considered.

In the following it is assumed that an initial iterate $x^0$ is
given and that iterates $x^k$ are generated via \mr{proc1}.
It is further assumed that the output of the algorithm
is given by weighted averages $\bar x^k$ defined as 
\begin{align}\label{barxk}
\bar x^k:= \left(\sum_{j=1}^k w_j\right)^{-1}\sum_{j=1}^k w_j x^j
\end{align}
and that
$\gamma_k\equiv \gamma(k)$ in \mr{proc1}, and $w_j \equiv w(j)$ where

\bee\label{specialcase}
\begin{array}{l}\gamma(t)\equiv c\left(\frac{M}{t+M}\right)^{\alpha},\\[4pt]
w(t) \equiv t^\beta,
\end{array}
\ \ \hbox{for some} \ \ \alpha\ge 0,\ \beta \ge 0,\ 0<c\le \frac 1{D_{n,n}}, \  M\ge 1.
\ene
The asymptotic analysis in \cite{polyak} (Assumption 2.2)
covers the case $\alpha\in(0,1)$, $\beta = 0$ and $M=1$.
Larger values of $\alpha$ lead to a faster reduction of the step length,
and for $M>1$ the step length reduction in the early iterations is slower.
In the numerical experiments the choice $M=1+\delta k^{max}$ is considered
where $k^{max}$ is the value of $k$ at which the iteration \mr{barxk}
is stopped, and $\delta \in[0,1]$ is a parameter to be determined.
If $\delta > 0$ is fixed, the step lengths $\gamma$ do not
converge to zero when $k^{max}\to\infty$.
\ignore{For $\delta = 0$, a choice $\alpha >1$ is not meaningful
since such step lengths are ``too short'' in the sense that
the iterates $x^k$ are not guaranteed to converge to the optimal
solution for $k^{max}\to\infty$ even in the ideal case of ``zero noise and 
$f(x)\equiv \frac 12 x^\T x$''.  
(The second part of Assumption 2.2 in \cite{polyak} is violated.)
For $\delta > 0$, however, $\gamma$ is
bounded below by some $\epsilon >0$ also for values of $\alpha >1$
so that the second part of Assumption 2.2 in \cite{polyak} holds.
In this paper an asymptotic analysis is not considered -- this aspect has been answered in \cite{polyak}. Instead, }
It is the aim of this paper
to determine optimal values $\alpha,\ c,\ M$ of the
step lengths and an optimal value $\beta$ of the weights
for finite values of $k^{max}$.

For completeness, the algorithm outlined above is
stated in detail:

\begin{algorithm}\label{alg1}\hfill\break
    Input:
    $x^0\in\re^n$, \ $\gamma(\, .\, ), \ w(\, .\, )$ as in \mr{specialcase},
    constant batch size $\nu$,
    number of iterations $k^{max}$.\hfill\break
    Initialization:
    Set $\sigma :=0$ (sum of the weights) and $\bar x^0:=0$
    (weighted sum of iterates).\hfill\break
    For $k=0,1,\ldots,k^{max}-1$ do \phantom{$\displaystyle\sum^K$}
    \begin{enumerate}
      \item Select a batch $S_k$ of size
        $|S_k|=\nu$ i.i.d.\ from $\{1,\ldots,m\}$.
      \item Set 
        $x^{k+1} := x^{k}-\gamma(k)\nabla f_{S_k}(x^k)$.
        \item Set $\bar x^{k+1} := \bar x^{k}+w(k+1) x^{k+1}$ and $\sigma := \sigma+w(k+1)$. 
    \end{enumerate}
    Set the weighted average $\bar x^{final}:=\bar x^{k^{max}}/\sigma $
    and return $\bar x^{final}$ as approximate minimizer of $f$.
\end{algorithm}
  
For $\alpha=\beta=0$, Algorithm \ref{alg1} reduces to an algorithm with
asymptotically optimal parameters in \cite{polyak}.
In order to identify optimal parameters depending on the number of iterations
allocated in advance
an elementary analysis of the iterates is presented next.

\section{Optimal Selection of Parameters}\label{sec.opt}

\subsection{Analysis Without Averaging}
With the above assumptions,
$$
x^1 =x^0-\gamma_0(\nabla f(x^0)+\noise^0) = x^0-\gamma_0Dx^0-\gamma_0\noise^0  =  (I-\gamma_{0} D) x^{0}  - \gamma_{0} \noise^{0}
$$
where $ (I-\gamma_{0} D)$ is a contraction satisfying
$$
0\preceq (1- \gamma_{0}D_{n,n}) I\preceq (I-\gamma_{0} D) \preceq
(1- \gamma_{0}D_{1,1}) I \prec I.
$$
For the next step, there is new noise ``$- \gamma_{1} \noise^{1}$''
while the noise
of the previous iteration is reduced,
\begin{align*}
x^2  &= (I-\gamma_{1} D) x^{1}  - \gamma_{1} \noise^{1}
= (I-\gamma_{1} D) [(I-\gamma_{0} D) x^{0}  - \gamma_{0} \noise^{0}]  - \gamma_{1} \noise^{1}\\
&=(I-\gamma_{1} D) (I-\gamma_{0} D) x^{0}  -
(I-\gamma_{1} D)\gamma_{0} \noise^{0}- \gamma_{1} \noise^{1}.
\end{align*}
Denote the contraction $I-\gamma_{i} D$ by $C_i$ and
observe that $\gamma_i\in(0,1/D_{n,n}]$ implies $||C_i||=1-\gamma_iD_{1,1} <1$.

Let the product of contractions be denoted by
\bee\label{Cik}
C_{i,k}=  \prod_{j=i}^{k} C_j=  \prod_{j=i}^{k} (I-\gamma_iD), 
\ene
where the empty product $C_{i,k}$ for $i>k$ is equal to $I$ by convention,
for example, $C_{k+1,k}=I$. Since the  largest entry of all $C_i$
  is at the first diagonal position, it follows that
$ ||C_{i,k}||=  \prod_{j=i}^{k} (1-\gamma_jD_{1,1})$.
Again, the empty product is 1, i.e.\ $\|C_{k+1,k}\|=1$.
The process \mr{proc1} then leads to the aggregated representation
\begin{align}\label{xk}
  x^k = \left[\prod_{i=0}^{k-1}(I-\gamma_iD)\right]x^0-\sum_{i=0}^{k-1}
  \left[\prod_{j=i+1}^{k-1}(I-\gamma_jD)\right]\gamma_i\noise^i
  = C_{0,k-1}x^0-\sum_{i=0}^{k-1}\gamma_{i} C_{i+1,k-1}\noise^{i}
\end{align}
with the expected value of $x^k$ given by
$E(x^k) = C_{0,k-1}x^0$
and noise term $\sum_{i=0}^{k-1}\gamma_{i} C_{i+1,k-1}\noise^{i}$.
  Observe that $E(x^k)$ coincides with the ``ideal'' iterates, meaning the iterates without noise
  and denote the ``ideal'' iterates by
  \bee\label{xkbb}
  \breve x^k := C_{0,k-1}x^0  \quad \hbox{for} \ \ k\ge 1.
  \ene

\subsection{Analysis With Averaging}
Now, consider weighted averaged iterates \mr{barxk}
for the weight function $w: \re_+\to\re_+$
given in \mr{specialcase} and $w_j=w(j)$.
Denote the ``ideal'' averages (without noise) by 
\begin{align}\label{xkb}
\breve{\bar x}^k := \left(\sum_{j=1}^k w_j\right)^{-1}\sum_{j=1}^k w_j \breve x^j
= \left(\sum_{j=1}^k w_j\right)^{-1}\sum_{j=1}^k w_j C_{0,j-1}x^0
\end{align}
and the ``accumulated'' noise in $x^k$ by
$$
\anoise^k := x^k-\breve x^k.
$$
The noise satisfies the recurrence relation $\anoise^0 :=0$ and
\begin{eqnarray}
\anoise^{k+1}
&=& x^{k+1}-\breve{x}^{k+1}
\overset{\eqref{xkbb}}{=} x^{k+1}-C_{0,k}x^0\nonumber \\
&\overset{\eqref{xk}}{=}& C_{0,k}x^0 -\sum_{i=0}^{k} \gamma_{i}C_{i+1,k}\,\noise^{i} - C_{0,k}x^0
= -\sum_{i=0}^{k} \gamma_{i}C_{i+1,k}\,\noise^{i}. \label{noisekpo}
\end{eqnarray}
For indices
$i\in\{0,\dots,k-1\}$ and $j\in\{1,\dots,k\}$ with $k\in\mathbb{N}$
a reordering of the sum 
 \begin{align}\label{switchsums}
  \sum_{j=1}^{k}\sum_{i=0}^{j-1}a_{i,j}=\sum_{i=0}^{k-1}\sum_{j=i+1}^{k}a_{i,j}.
 \end{align}
 leads to a representation of the ``total''
 noise $\bar x^k-\breve{\bar x}^k$ in $\bar x^k$  denoted by
$\tnoise^k$,
\begin{eqnarray*}
  \tnoise^k:=\bar x^k- \breve{\bar x}^k
  &=& \left(\sum_{j=1}^k w_j\right)^{-1}\sum_{j=1}^k w_j (x^j-\breve x^j)\\
  & \overset{\eqref{noisekpo}}{=}  & -\left(\sum_{j=1}^k w_j\right)^{-1}\sum_{j=1}^k w_j
  \sum_{i=0}^{j-1} \gamma_{i}C_{i+1,j-1}\,\noise^{i}\\
  &\overset{\eqref{switchsums}}{=}& -\left(\sum_{j=1}^k w_j\right)^{-1}\sum_{i=0}^{k-1}
  \sum_{j=i+1}^k w_j \gamma_{i}C_{i+1,j-1}\, \noise^{i}\\
 &=&\left(\sum_{j=1}^k w_j\right)^{-1}\sum_{i=0}^{k-1}
   G_{i,k} \, \noise^{i}\\ 
\end{eqnarray*}
where
\bee \label{Mik}
G_{i,k}:=-\gamma_{i}\sum_{j=i+1}^k w_j C_{i+1,j-1}
\qquad \hbox{for} \qquad 0\le i\le k-1.
\ene

\subsection{Parameter Selection Without Averaging}
For a comparison, first an algorithm is considered that generates the ``plain''
iterate $x^{k^{max}}$ as output rather than the weighted average $\bar x^{final}$.

In this case, based on \mr{xk}, the quantity to be minimized
in the design of a method with optimal convergence would be a weighted sum of $\|C_{0,k-1}\|$ and
of a bound of the variance of the noise term.
Due to stochastic independence of the noise terms  \mr{iid},
  the covariance matrix of $x^k$ satisfies
$$
\text{Cov}\left(\sum_{i=0}^{k-1}\gamma_{i} C_{i+1,k-1}\noise^{i}\right)
= \sum_{i=0}^{k-1}\gamma_{i}^2
\left(C_{i+1,k-1}\right)\text{Cov}(\noise^{i})\left(C_{i+1,k-1}\right)^\text{T}.
$$
Furthermore, by \mr{iid} the covariance matrices $\text{Cov}(\noise^{i})$ are all equal,
$\text{Cov}(\noise^{i})\equiv \Sigma$ for some positive semidefinite $\Sigma$.
In the case that $\Sigma$ is a multiple of the identity matrix, one obtains an exact reformulation
$$ 
   \left\| \text{Cov}\left(\sum_{i=0}^{k-1}\gamma_{i} C_{i+1,k-1}\noise^{i}\right)\right\|
   = \left\| \sum_{i=0}^{k-1}\gamma_{i}^2
\left(C_{i+1,k-1}\right)\Sigma\left(C_{i+1,k-1}\right) \right\|
= \sum_{i=0}^{k-1}\gamma_{i}^2 \|C_{i+1,k-1}\|^2\|\Sigma\|
$$
since all $C_{i+1,k-1}$ are diagonal matrices with their largest entry
defining their norm at the $(1,1)$-position,
and in the case of a general positive semidefinite matrix $\Sigma$ one obtains an upper bound
\bee\label{cov}  
   \left\| \text{Cov}\left(\sum_{i=0}^{k-1}\gamma_{i} C_{i+1,k-1}\noise^{i}\right)\right\|
   = \left\| \sum_{i=0}^{k-1}\gamma_{i}^2
\left(C_{i+1,k-1}\right)\Sigma\left(C_{i+1,k-1}\right) \right\|
\leq \sum_{i=0}^{k-1}\gamma_{i}^2 \|C_{i+1,k-1}\|^2\|\Sigma\|
\ene
  where the inequality follows from a repeated application of the triangle
  inequality and from sub-multiplicativity of the norm.
  Let
\begin{align}\label{theta}
\vartheta(\gamma) := \left(\sum_{i=0}^{k-1}\gamma_{i}^2 \left\|C_{i+1,k-1}\right\|^2\right)^{1/2}.
\end{align}
Optimizing the
expected norm of $x^k$
\bee\label{explainmu}
E\left(\norm{x^k}\right)\le \| C_{0,k-1}\| \|x^0\| +\vartheta(\gamma)\left\| \Sigma \right\|^{1/2}
\ene
thus leads to the aim of choosing the step length function $\gamma(\, .\, )$ such that
\bee\label{thetamin}
\left(\,\|C_{0,k-1}\|+\mu \, \vartheta(\gamma)\,\right)/(1+\mu)
\ene
is minimized   
where $\mu\ge 0$ is a fixed weight.
For $\mu:= \left\| \Sigma \right\|^{1/2}/\|x^0\|$
the minimizers of \mr{explainmu} and \mr{thetamin} coincide; but unfortunately, the ratio of
``noise'' $\left\| \Sigma \right\|^{1/2}$ to ``starting error''
$\|x^0\|$ generally is not known. Nevertheless, the separation
of ``optimization error'' (here $\| C_{0,k-1}\| \|x^0\|$) and ``stochastic error''
(here $\vartheta(\gamma)\left\| \Sigma \right\|^{1/2}$)
can be extended to weighted averages in the next subsection and will then be exploited
with the aim of identifying suitable parameters $\alpha,\beta,\delta,M$.
\ignore{The choice of $\mu$, however, is not evident.
For large ``initial errors $\|x^0\|$'' smaller values of $\mu$ appear to be appropriate
(to concentrate on the minimization of the optimization error) and for
large noise terms (i.e.\ large $\|\Sigma\|^{1/2}$) larger values of $\mu$ appear to be
appropriate. Unfortunately the ratio $\|\Sigma\|^{1/2}/\|x^0\|$ is not known in general.}

\subsection{Parameter Selection With Averaging}
Since the step lengths and weights are pre-defined the same is true for
$G_{i,k}$, and thus, since the noise terms $\noise^{i}$ are assumed
to be i.i.d., also $G_{i,k}\noise^{i}$
are independently distributed.
The variance of $\tnoise^k$ therefore is $\left(\sum_{j=1}^k w_j\right)^{-2}$
times the sum of the variances of $G_{i,k} \, \noise^{i}$.
The latter are bounded by a fixed multiple (depending on the distribution of
the terms $\noise^{i}$) of $\|G_{i,k}\|^2$.

To reduce the expected value of $||\tnoise^k||_2$, it is therefore the aim to
define the weight function $w$ and the step length function $\gamma$
such that 
\begin{align}\label{kappa}
\kappa(w,\gamma):=
\left(\left(\sum_{j=1}^k w_j\right)^{\!\!-2}\,\sum_{i=0}^{k-1}   ||G_{i,k}||^{2}
\right)^{1/2}
\end{align}
is small. As in the case of $\vartheta$ in \mr{theta}, the
  upper bound $\kappa^2$ is
  the exact norm of the covariance matrix of $\tnoise^k$ when $\Sigma$ is a
  multiple of the identity matrix.

Simultaneously, the second goal is that also the norm of
$\breve{\bar x}^k$ should be small, i.e.
\begin{align}\label{tau}
\tau(w,\gamma):=
\left\|
\left(\sum_{j=1}^k w_j\right)^{\!\!-1}\,\sum_{j=1}^k w_j
\prod_{i=0}^{j-1}(I-\gamma_i D)
\right\|
=
\left(\sum_{j=1}^k w_j\right)^{\!\!-1}\,\sum_{j=1}^k w_j
\left\|C_{0,j-1}
\right\|
\end{align}
should be small.
The above equation again uses the fact that all $C_i$ are diagonal
matrices with the largest entry defining the norm at the $(1,1)$ position.
Thus, in place of \mr{thetamin} it is the aim of choosing the step length
function $\gamma(\, .\, )$ and the weight function $w(\, .\, )$ such that
\bee\label{obj.fn}
r(w,\gamma):= \left(\,\tau(w,\gamma) + \mu \kappa(w,\gamma)\,\right)/(1+\mu)
\ene
is minimized for a given fixed $\mu\ge 0$.
Again, an ``appropriate'' choice of $\mu$ is not evident.

For the special case \mr{specialcase}
the evaluation of $\kappa$ and $\tau$ with order
$k$ arithmetic operations is considered next.

\subsubsection{Evaluation of All Norms  $\|C_{0,j}\|$ and $\|C_{i+1,k-1}\|$}
The choice of $\gamma\in (0, 1/D_{n,n}]$ implies that
\begin{align}\label{normIminusGammaD}
  \| C_k \|_2 =
  \| I-\gamma_k D\|_2 = 1-\gamma_k D_{1,1} = 1- \bar c \left(\frac{M}{k+M}\right)^{\!\alpha} \quad
\hbox{for all} \ k\ge 0
\end{align}
where
$$
\bar c := cD_{1,1} \le \frac{D_{1,1}}{D_{n,n}} = \frac 1 {\hbox{cond}(D)}.
$$
For  $0\le j \le k$ note that 
$$
\|C_{0,j}\|=\|\prod_{\ell = 0}^j C_\ell \| 
=\prod_{\ell = 0}^j \!\left(1-\bar c\left(\frac{M}{\ell+M}\right)^{\!\alpha}\right)
=\|C_{0,j-1}\|\!\left(1-\bar c\left(\frac{M}{j+M}\right)^{\!\alpha}   \right)
$$ 
where the second equality again follows from the diagonal
structure of $C_\ell=I-\gamma_\ell D$ with the largest entry
in absolute value always at the (1,1)-position.

Thus, starting from $\|C_{0,0}\|=1-\bar c$,
all $\|C_{0,j}\|$ can be computed for $j=1,2,3\ldots,k$ with a total of
order $k$ arithmetic operations.\\
Likewise, starting with $\|C_{k,k-1}\|$, which is 1 by convention, the predecessors
  $\|C_{j,k-1}\|$ can be iteratively determined by
\begin{align*}
\|C_{j,k-1}\| = \prod_{\ell = j}^{k-1} \left(1-\bar c\left(\frac{M}{\ell+M}\right)^{\!\alpha}   \right)
= \|C_{j+1,k-1}\| \left(1-\bar c\left(\frac{M}{j+M}\right)^{\!\alpha} \right)
\end{align*}
for $j\in\{k-1,\dots,0\}$, and $\vartheta$ in \mr{theta} can be evaluated
with order $k$ arithmetic operations.

\subsubsection{Evaluation of All Norms $\|G_{i,k}\|$}
The aim of this section is to derive a scheme
for evaluating all norms $\|G_{i,k}\|$ for $0\le i\le k-1$
also with order $k$ arithmetic operations, so that
the parameters $\alpha$, $\beta$, $c$, and $M$ can easily be optimized for
maximum iteration numbers $k^{max}$ up to the order of about $10^8$.

Note that for $0\le i\le k-1$:
\bee\label{miknew}
\|G_{i,k}\| = \left\| \gamma_i\!\!\sum_{j=i+1}^{k}w_jC_{i+1,j-1} \right\| 
  = \gamma_i\left\| \sum_{j=i+1}^{k}j^\beta C_{i+1,j-1} \right\| 
  = c \left(\frac{M}  {i+M}\right)^{\!\alpha} \!\!\sum_{j=i+1}^{k}j^\beta \| C_{i+1,j-1} \|
\ene
where the last equality follows again since all $C_i$ are diagonal matrices
with the largest entry defining the norm at the (1,1) position.
Since $\| C_{i+1,i} \|=1$ it follows that
$$
\sum_{j=i+1}^{k}j^\beta \| C_{i+1,j-1} \|
=
(i+1)^\beta+\sum_{j=i+2}^{k}j^\beta \| C_{i+1,j-1} \|
=
(i+1)^\beta+\sum_{j=i+2}^{k}j^\beta \|C_{i+1}\|\, \| C_{i+2,j-1} \|
$$ 
$$
=
(i+1)^\beta+\|C_{i+1}\|\sum_{j=i+2}^{k}j^\beta  \| C_{i+2,j-1} \|
=(i+1)^\beta+\|C_{i+1}\|\left(\frac{(i+1+M)^\alpha}{cM^\alpha} \|G_{i+1,k}\|\right).
$$
In the last equation relation \mr{miknew} has been used for
$\|G_{i+1,k}\|$ in place of $\|G_{i,k}\|$.
Hence, starting with
$\|G_{k-1,k}\|=\frac{ck^{\beta}M^\alpha}{(k-1+M)^\alpha}$,
all $\|G_{i,k}\|$
can be computed for $i=k-2$, $k-3$, \ldots
via
$$
\|G_{i,k}\| =c\left(\frac{M}{i+M}\right)^\alpha \!\! (i+1)^{\beta}\, +\, \|C_{i+1}\|
\left(\frac{i+1+M}{i+M}\right)^\alpha \|G_{i+1,k}\|
$$
with $\|C_{i+1}\|=1-\frac{\bar cM^\alpha}{(i+1+M)^\alpha}$.

For the case $\alpha = \beta = 0$ as considered in \cite{polyak},
  the above simplifies to \\
  $\|G_{k-j,k}\|=c(1-(1-\bar c)^j)/\bar c$, and
  the quantities $\kappa$ and $\tau$ allow a closed form representation,
  $$
  \kappa = \frac {c}{k^{max}\,\bar c}\left(k^{max}-
  2\frac{1-\bar c-(1-\bar c)^{k^{max}+1}}{\bar c}+
  \frac{(1-\bar c)^2-(1-\bar c)^{2k^{max}+2}}{1-(1-\bar c)^2}
  \right)^{1/2}
  $$
  and
  $$
\tau = \frac{(1-\bar c)(1-(1-\bar c)^{k^{max}})}{k^{max}\, \bar c}.
$$
The straightforward derivation of these formulas is omitted for brevity.
For fixed values of $\bar c>0$ it follows that $\kappa$ is of the order
$1/\sqrt{k^{max}}$ and $\tau$ is of the order $1/k^{max}$ so that
convergence of Algorithm \ref{alg1} to the
optimal solution when $k^{max}\to\infty$
follows also for the case $\alpha=\beta=0$.

With these preparations it is now possible to evaluate
$\tau(w,\gamma)$ and $\kappa(w,\gamma)$ in \mr{tau} and \mr{kappa}
with order $k$ arithmetic operations and thus to minimize the function $r$ with 
$w,\gamma$ as in \mr{specialcase}.
Here, larger values of $\mu$ are meaningful, when the noise is
large compared to the distance of the initial iterate from optimality.

Note that when $f$ is multiplied by some constant $\eta>0$
then $\bar c$ and all $\|C_{i,j}\|$ remain invariant, but 
$c$ and all $\|G_{i,k}\|$ are multiplied by $1/\eta$.
Hence, when multiplying $\mu$ with $\eta$, the minimizer of
\mr{obj.fn} remains invariant.

\vfill\eject

\section{Numerical Examples} \label{sec.num}

\subsection{Optimizing the Parameters of Weighted Averaging}
\label{sec.num1}

In this sub-section the selection of the parameters $\alpha,\ \beta,\ c$, and $M$
in algorithm \ref{alg1} is considered minimizing the function $r$ in \mr{obj.fn}.
Here, $r$ is a weighted sum of $\tau$ and $\kappa$
where $\tau\ge 0$ always is less than 1 and decreases with increasing
values of $\beta$.
On the other hand, $\kappa$ is minimized for $\beta = 0$.
It is the aim of the considerations below to balance the conflicting goals of
minimizing both $\tau$ and $\kappa$. 

To standardize the results in this sub-section the largest eigenvalue
  of $f$ is fixed to
  $$
D_{n,n}=1
  $$
throughout.

\bigskip

Table~1 and Table~2 give some intuition about the values of 
$\tau$ and $\kappa$ when $\alpha =\beta=0$ and $c$ is also fixed to
$c=1$.

\bigskip
\vbox{
\begin{center}
\footnotesize{
\begin{tabular}{|c|c|c|c|c|c|c|c|c|}
\hline  
\diagbox{$k^{max}$}{$cond(D)$}&  $10^{0.5}$ & $10^{1}$ & $10^{1.5}$ & $10^{2}$ &  $10^{2.5}$ & $10^{3}$ & $10^{3.5}$ & $10^{4}$\\
\hline 
$10^{2}$         &      -1.6646  &   -1.0458  &   -0.5315  &   -0.2023  &   -0.0676  &   -0.0218  &   -0.0069  &   -0.0022\\
$10^{2.5}$        &      -2.1643  &   -1.5454  &   -1.0133  &   -0.5226  &   -0.1999  &   -0.0671  &   -0.0216  &   -0.0069\\
$10^{3}$         &      -2.6646  &   -2.0458  &   -1.5136  &   -1.0044  &   -0.5198  &   -0.1995  &   -0.0669  &   -0.0216\\
$10^{3.5}$        &      -3.1643  &   -2.5454  &   -2.0133  &   -1.5041  &   -1.0008  &   -0.5189  &   -0.1990  &   -0.0668\\
$10^{4}$         &      -3.6646  &   -3.0458  &   -2.5136  &   -2.0044  &   -1.5011  &   -1.0005  &   -0.5186  &   -0.1992\\
$10^{4.5}$        &      -4.1643  &   -3.5454  &   -3.0133  &   -2.5041  &   -2.0007  &   -1.5001  &   -0.9995  &   -0.5186\\
$10^{5}$         &      -4.6646  &   -4.0458  &   -3.5136  &   -3.0044  &   -2.5011  &   -2.0004  &   -1.4998  &   -1.0001\\
$10^{5.5}$        &      -5.1643  &   -4.5454  &   -4.0133  &   -3.5041  &   -3.0007  &   -2.5001  &   -1.9995  &   -1.4997\\
$10^{6}$         &      -5.6646  &   -5.0458  &   -4.5136  &   -4.0044  &   -3.5011  &   -3.0004  &   -2.4998  &   -2.0000\\
$10^{6.5}$        &      -6.1643  &   -5.5454  &   -5.0133  &   -4.5041  &   -4.0007  &   -3.5001  &   -2.9995  &   -2.4997\\
$10^{7}$         &      -6.6646  &   -6.0458  &   -5.5136  &   -5.0044  &   -4.5011  &   -4.0004  &   -3.4998  &   -3.0000\\
$10^{7.5}$        &      -7.1643  &   -6.5454  &   -6.0133  &   -5.5041  &   -5.0007  &   -4.5001  &   -3.9995  &   -3.4997\\
$10^{8}$         &      -7.6646  &   -7.0458  &   -6.5136  &   -6.0044  &   -5.5011  &   -5.0004  &   -4.4998  &   -4.0000\\
\hline
\end{tabular}
}
\vskip 12pt
\centerline{{\bf Table~1:} Values of $\log_{10}(\tau(0,0))$ for
    $k^{max} =  10^2,\ 10^{2.5 },\ 10^3,\ \ldots , \ 10^8$ in rows 1 -- 13}
\centerline{and condition numbers $10^{0.5},\ 10^1,\ 10^{1.5}, \  \ldots 10^4$ in columns 1 -- 7}
  \end{center}
}

\bigskip

\noindent The values of $\tau(0,0)$ in Table~1 are not surprising,
  and are listed only as a comparison to Table~2 below.
  When $k^{max}$ is larger than the condition number then $\tau(0,0)$ is
  roughly of the order $\frac{\hbox{condition number}}{k^{max}}$.

\bigskip

In Table~2 below
  it is interesting to observe that for large condition numbers 
  such as condition number $10^4$ in column 7, the bound $\kappa(0,0)$
  for the variance increases first when  $k^{max}$ increases, and
  starting from $k^{max} = 100$ it reaches a maximum of
  $\kappa(0,0) \approx 41.7920$ for $k^{max} = 10^{4.5}$ in row 6
  before decreasing again for larger values of $k^{max}$.
  (For the dimension, and below also for $k^{max}$,
  square roots are rounded to integer values
  in Table~1 and Table~2.)

\vbox{
\begin{center}
\footnotesize{
\begin{tabular}{|c|c|c|c|c|c|c|c|c|}
\hline 
\diagbox{$k^{max}$}{$cond(D)$} &  $10^{0.5}$ & $10^{1}$ & $10^{1.5}$ & $10^{2}$ &  $10^{2.5}$ & $10^{3}$ & $10^{3.5}$ & $10^{4}$\\
\hline 
$10^{2}$   &       0.3109 &   0.9288  &  2.3730 &   4.1385 &   5.1887 &   5.6063 &   5.7490 &   5.7952\\
$10^{2.5 }$ &       0.1770 &   0.5502 &   1.6449 &   4.1915 &   7.3131 &   9.1714 &   9.9138 &   10.1670\\
$10^{3}$   &       0.0999 &   0.3140 &   0.9773 &   2.9176 &   7.4412 &   12.9772 &   16.2857 &   17.6052\\
$10^{3.5 }$ &       0.0563 &   0.1775 &   0.5588 &   1.7365 &   5.1907 &   13.2193 &   23.0643 &   28.9336\\
$10^{4}$   &       0.0316 &   0.0999 &   0.3157 &   0.9925 &   3.0887 &   9.2203 &   23.5113 &   41.0028\\
$10^{4.5 }$ &       0.0178 &   0.0562 &   0.1779 &   0.5612 &   1.7668 &   5.4904 &   16.4110 &   41.7920\\
$10^{5}$   &       0.0100 &   0.0316 &   0.1000 &   0.3160 &   0.9983 &   3.1385 &   9.7669 &   29.1551\\
$10^{5.5 }$ &       0.0056 &   0.0178 &   0.0563 &   0.1779 &   0.5625 &   1.7747 &   5.5871 &   17.3619\\
$10^{6}$   &       0.0032 &   0.0100 &   0.0316 &   0.1000 &   0.3164 &   0.9993 &   3.1570 &   9.9247\\
$10^{6.5 }$ &       0.0018 &   0.0056 &   0.0178 &   0.0563 &   0.1780 &   0.5624 &   1.7789 &   5.6121\\
$10^{7}$   &       0.0010 &   0.0032 &   0.0100 &   0.0316 &   0.1001 &   0.3162 &   1.0005 &   3.1599\\
$10^{7.5 }$ &       0.0006 &   0.0018 &   0.0056 &   0.0178 &   0.0563 &   0.1779 &   0.5629 &   1.7785\\
$10^{8}$   &       0.0003 &   0.0010 &   0.0032 &   0.0100 &   0.0316 &   0.1000 &   0.3164 &   0.9999\\
\hline
\end{tabular}
}
\vskip 12pt
\centerline{{\bf Table~2:} Values of $\kappa(0,0)$ for
    $k^{max} =  10^2,\ 10^{2.5 },\ 10^3,\ \ldots , \ 10^8$ in rows 1 -- 13}
\centerline{and condition numbers $10^{0.5},\ 10^1,\ 10^{1.5}, \  \ldots 10^4$ in columns 1 -- 7}
  \end{center}
}

\bigskip

While condition numbers as large as $10^4$ might not be typical for stochastic applications,
and while moderate iteration numbers such as $k^{max}=10^{4.5}$ cannot render significant progress for such large
condition numbers, 
it is interesting to consider possible improvements of $\kappa$ for $k^{max} = 10^{4.5}$ and $D_{1,1}=10^{-4}$ by
optimizing $\tau + \mu\, \kappa$ with respect to $\alpha$, $\beta$, $c$, and $M$ for different values of $\mu$.

For minimizing the function $r=(\tau+\mu\, \kappa)/(1+\mu)$ in \mr{obj.fn}
the descent algorithm
``min$\underline{\ }$f.m'' from \cite{lazar}  was used that
aims for a local minimizer near the starting point:

The variables $\alpha$ and $\beta$ were constrained to the intervals
$[0,\ 2]$ and $[0,\ 5]$. (The upper bounds $\alpha \le 2$ and  $\beta \le 5$
were chosen at will to limit the search space to a compact domain.)
For $c\le 1$ a lower bound of 0.1 was chosen.
(The reduction of $\tau$ is considered as too slow when $c<0.1$.)
Finally $M$ was set as $1+\delta k^{max}$ with $\delta \in [0,1]$.
In Table~3, the (approximately) optimal parameters identified
with min$\underline{\ }$f
are listed for different weights $\mu > 0$ and $k^{max} = 10^{4.5}$, $D_{1,1}=10^{-4}$.

For each run of min$\underline{\ }$f, the four starting values
$(\alpha,\beta) = (0,0),\ (0,2),\ (\frac 12,0),\ (\frac 12,2)$
were used as well as $\delta = 0.1$ and $c=0.5$.
When min$\underline{\ }$f identified different
approximate optimal solutions the one with the lowest value of
$r(\alpha,\beta)$ is listed.

\vbox{
\begin{center}
\begin{tabular}{|c|c|c|c|c|c|c|c|}
\hline
$\mu$       &   NA      &  1   & 0.1     & 0.017 & 0.017 & 0.01       &  0.001   \\
\hline
$\alpha$    &   0       &  2     & 2      & 2     & 0     &  0         &  0       \\
$\beta$     &   0       &  0     & 0      & 0     & 0.718 &  2.081     &  5       \\
$c$         &   1       &  0.1   & 0.1    & 0.1   & 1     &  1         &  1       \\
$\delta$    &   0       &  0     & 0      & 0     & 0     &  0         &  0       \\
$\tau$      &   0.303   &  1.000 & 1.000  & 1.000 & 0.189 &  0.114     &  0.073   \\
$\kappa$    &   41.79   &  0.104 & 0.104  & 0.104 & 47.60 &  53.18     &  58.84   \\
$r$         &           &  0.552 & 0.919  & 0.985 & 0.982 &  0.639     &  0.132   \\
\hline
\end{tabular}
\vskip 12pt
\centerline{{\bf Table 3:} Optimal parameters $\alpha$, $\beta$, $c$,
  $\delta$ for $k^{max} = 10^{4.5}$, $D_{1,1}=10^{-4}$, and different values of $\mu$.}
\end{center}
}

The numbers in the first column refer to the situation of Polyak and Juditsky \cite{polyak}
and are not the result of an optimization process. The two columns for $\mu = 0.017$ show the
situation that two approximate local minimizers were found with similar objective value ``$r$''
but rather different input arguments.
The optimal parameters appear to be discontinuous near $\mu = 0.017$.
Since the ``appropriate'' weight $\mu> 0$ depends on the unknown distance of $x^0$ to the
optimal solution and on the unknown magnitude of the noise, the selection of
$\alpha$, $\beta$, $c$, and $\delta$ cannot be extracted from
data such as Table~3, even if $D_{1,1}$ and $D_{n,n}$ are known.

Given an example that is not as poorly conditioned as in Table~3, namely
$$
D_{1,1} = 0.03,
$$
different values of $k^{max}$ and weighting terms $\mu$ were considered
for Table~4 such that 
at least two approximate local optimal solutions $\alpha$, $\beta$, $c$, $\delta$
could be identified. 

\vbox{
\begin{center}
\begin{tabular}{|c||c|c||c|c||c|c|c|}
  \hline
$k^{max}$    &   \multicolumn{2}{c||}{1000}    &   \multicolumn{2}{c||}{10000}    &  \multicolumn{3}{c|}{100000}  \\
$\mu$       &   \multicolumn{2}{c||}{0.05}    &   \multicolumn{2}{c||}{0.012}    &  \multicolumn{3}{c|}{0.00148} \\
\hline
$\alpha$    &   1.104   &  0       &   0.519    &  0        &  1.195   & 0        &  0       \\
$\beta$     &   1.382   &  0.809   &   0.614    &  0.606    &  5       & 0.5955   &  0.6521  \\
$c$         &   1       &  1       &   1        &  1        &  1       & 1        &  1       \\
$\delta$    &   0.186   &  0       &   0.164    &  0        &  1.29e-3 & 0        &  0       \\
$\tau$      &   2.37e-3 &  3.49e-3 &   0.3594   &  0.3593   &  5.92e-7 & 3.94e-6  &  2.61e-6 \\
$\kappa$    &   1.171   &  1.152   &   1.45e-4  &  1.47e-4  &  0.119   & 0.114    &  0.115   \\
$r$         &   5.80e-2 &  5.82e-2 &   4.4054e-3&  4.4057e-3&  1.66e-4 & 1.72e-4  &  1.72e-4 \\
\hline
\end{tabular}
\vskip 12pt
\centerline{{\bf Table 4:} Nearly optimal parameters $\alpha$, $\beta$, $c$,
  $\delta$ for  different values of $\mu$ and $k^{max}$.}
\end{center}
}

\bigskip

Because of the discontinuous dependence of $\alpha$ and $\beta$ (and also of $\tau$ and $\kappa$)
on $\mu$ as observed in Table~3 and Table~4, another selection process
for optimizing the parameters $\alpha,\beta,c,\delta$ was considered, namely 
\bee\label{nlp}
\min_{-1\le v_1,v_2 \le 0.1}  v_1 +\mu v_2 \ | \  \tau(w,\gamma) = (1+v_1)\tau(w^0,\gamma^0),\ \
\kappa(w,\gamma) = (1+v_2)\kappa(w^0,\gamma^0).
\ene
Problem \mr{nlp} uses a compact notation highlighting the changes compared to \mr {obj.fn}
in Table~3 and Table~4.
As in Table~3 and Table~4, $w$ depends on $\beta\in[0,5]$ and $\gamma$
depends on $\alpha\in[0,2]$,
$c\in[0.1,1]$, and $\delta\in[0,1]$,
and $w^0$, $\gamma^0$ refer to $\alpha = \beta = \delta = 0$, and  $c=1$. 
The restriction $v_1,v_2\le 0.1$ implies that neither $\tau$ nor $\kappa$ are allowed
to be more than 10\% worse than the choice $w^0,\gamma^0$.
Again, $\mu>0$ is a weight balancing optimization error and stochastic error;
however, due to the upper bound on $v_1,v_2$ the dependence on $\mu$ of the optimal solution
turns out to be less pronounced.

Using ``min$\underline{\ }$fc.m'' from \cite{jarrelieder}
with different starting points, approximate local optimal solutions
were computed for Problem \mr{nlp}. For upper bounds $v_1,v_2\le 0$
aiming at a simultaneous reduction of both $\tau$ and $\kappa$, no point apart from $\alpha=\beta=0$
could be identified via ``min$\underline{\ }$fc.m''. However, allowing a small slack of
$v_1,v_2\le 0.1$
as in \mr{nlp} led to
\bee\label{parset}
\alpha = \delta = 0, \ \  \beta = 0.7116,\ \  \hbox{and} \ \  c=1
\ene
for  a ``generic''
situation with $D_{1,1}=0.03$, $D_{n,n}=1$ and $k^{max}=10000$.
For these values of $\alpha,\beta,c,\delta$, the value of
$\kappa$ did increase by 10\% (from 0.3325 to 0.3658) while $\tau$ decreased by 97.4\%
(from 3.2e-3 to 8.7e-5) compared to the choice $\alpha=\beta=\delta=0$
and $C=1$.
The large value of $\kappa$ indicates for this condition number
that $10^4$ sgd-iterations
with noise do not reduce the expected error below approximately $1/3$
of the size of the noise-terms.
If $k^{max}$ is increased to $10^6$, the value of $\kappa$ reduces to
0.0366 and $\tau$ reduces to $3.3\cdot 10^{-8}$. 

For large initial errors $\|x^0\|$ this reduction of $\tau$
will outweigh the 10\%-increase of $\kappa$.
In this situation a significant gain in the optimization error could be
achieved when
allowing a small increase of the stochastic error.

To test the robustness of this solution, the other 
parameters $k^{max}$ and $D_{1,1}$ of Table~1 and Table~2 were evaluated as well
for $\alpha = \delta = 0$, $\beta = 0.7116$:
For higher condition numbers up to $10^4$ the
increase of $\kappa$ is less than 21\% and it is less than 10\% for smaller condition numbers
or larger values of $k^{max}$. Likewise, the reduction of $\tau$ deteriorates for increasing condition
numbers but improves for smaller condition numbers or for larger values of $k^{max}$.

\bigskip

\subsubsection{Practical Parameters}
Summarizing, while a step length reduction generally does not lead to a significant improvement
of convergence, moderately growing weights such as $w_j = j^{0.7}$ do lead to faster
reduction of the optimization error without deteriorating the stochastic error.
Thus, the parameter setting \mr{parset}
is used in the numerical examples in Subsection \ref{sec.test}
with the slight modification that $\beta$ is set to $\beta = 0.7$ (for simplicity).

\subsection{Test Examples}\label{sec.test}
Algorithm \ref{alg1} is rather simple and the results in
Section \ref{sec.opt} are
not very strong but ``robust'' in the sense that they
are independent of the dimension $n$ and of the number $m$ in the definition
\mr{f} of $f$.

\bigskip

\subsubsection{Simple Examples}
For the first set of examples a situation is considered that 
does satisfy the strong assumptions of Section \ref{sec.opt}
and where $m=\infty$ so that variance reduction by periodic
full gradient evaluations is not possible.

For a given dimension $n\in[10^2,10^8]$ a non-singular
randomly generated diagonal matrix $D\in\re^{n\times n}$ is chosen
with diagonal entries in $[0.1,1]$
and a scaling factor $\rho := 1/\sqrt n$.
For $k=1,2,\ldots$ random  vectors $b^k$ are drawn independently
from an $n$-variate normal
distribution with expected value $0_n$ and covariance matrix $\rho^2I_n$, i.\,e.\ 
$b^k\sim N_n(0_n,\rho^2I_n)$. Thus, $E(\|b^k\|_2^2)=n\rho^2=1$ independently of
the dimension $n$. The functions $f_k$ are then given by
$$
f_k(x)\equiv \tfrac 12 x^\T Dx+(b^k)^\T x
$$
with $E(f_k(x)) = \tfrac 12 x^\T Dx$.

For $k^{max}=10^5$ the results of Algorithm \ref{alg1} with
$\alpha = 0$, $\beta = 0.7$ and $x^0=e/\sqrt n$
(i.e. $\|x^0\|_2=1$) are listed for different dimensions:

\vbox{
\begin{center}
\begin{tabular}{|l||c|c|c|c|c|}
  \hline
$n$           &  $10^1$ & $10^2$ & $10^3$ & $10^4$ & $10^5$  \\
  \hline
$\|\bar x^{final}\|_2$ & 0.0123 & 0.0101 & 0.0129 & 0.0135 & 0.0133 \\
\hline
\end{tabular}
\vskip 12pt
\centerline{{\bf Table 5:} Final error in dependence of the dimension.}
\end{center}
}
For the examples listed above the number of iterations indeed does
not display any dependence on the dimension.

\bigskip

For the same setting with $n=100$ and different values of $k^{max}$
the results of Algorithm~\ref{alg1} are as follows:

\vbox{
\begin{center}
\begin{tabular}{|l||c|c|c|c|c|}
  \hline
$k^{max}$            & $10^4$ & $10^5$ & $10^6$ & $10^7$ & $10^8$    \\
  \hline
$\|\bar x^{final}\|_2$ & 3.9e-2 & 1.4e-2 & 3.8e-3 & 1.2e-3 & 3.9e-4\\
\hline
\end{tabular}
\vskip 12pt
\centerline{{\bf Table 6:} Final error in dependence of $k^{max}$.}
\end{center}
}
Here, the convergence is rather slow with a growing number
of iterations -- it is of the order $1/\sqrt{k^{max}}$, i.e., as
in the analysis of Polyak and Juditsky \cite{polyak},
the stochastic effects do dominate the convergence.
The reductions of the initial error
observed in Table~6 are better by a factor of 2.5 or 2.6 than the
theoretical bounds $\kappa(0,0)$ listed in the second column of Table~2.
  
For all test runs listed in Table~5 and Table~6,
the final iterate had a norm $\|x^{k^{max}}\|\in[1.21,1.39]$ ---
since $\alpha = 0$, the step length was constant and the
final iterate was largely determined by the noise.

\bigskip

The same setting with $n=100$ and $k^{max}=10^5$ was now
applied to different starting points $x^0 = \lambda e/\sqrt n$
with $\lambda >0$ and with $\beta = 0.7$
as well as $\beta = 0$.

\vbox{
\begin{center}
\begin{tabular}{|c|l||c|c|c|c|c|}
  \hline
\multicolumn{2}{|c||}{$\lambda$\ \ (i.e. $\|x^0\|_2$ )}  & $10^0$ & $10^2$ & $10^4$ & $10^6$ & $10^8$    \\
  \hline
\multirow{2}{*}{$\|\bar x^{final}\|_2$} & $\beta=0$   & 1.0e-2 & 1.1e-2 & 2.9e-1 & 28.6 & 3179\\
 & $\beta=0.7$ & 1.4e-2 & 1.2e-2 & 1.3e-2 & 5.4e-2 & 6.1\\
\hline
\end{tabular}
\vskip 12pt
\centerline{{\bf Table 7:} Final error in dependence of $\|x^0\|_2$.}
\end{center}
}
For a large initial error $\lambda$
there is a significant gain when choosing
$\beta = 0.7$ compared to $\beta=0$
while there is not much loss of the choice $\beta = 0.7$ for small
initial errors.
(For $\lambda=10^8$
the final iterate $x^{100000}$ has norm 1.26, i.e.\ it is much closer to the
optimal solution than the average $\|\bar x^{final}\|_2$ with $\beta=0$.)

For large initial errors, the value of $\tau$ will dominate
the convergence behavior. Indeed, 
for the entry with $\beta = 0.7$ in the last column of Table~7 a reduction
of the initial error by $10^8/6.1\approx 1.6\cdot 10^7$ can be observed
which is larger by a (moderate) factor of about 3.5 than the
theoretical bound $1/\tau(0,0.7) \approx 4.5\cdot 10^6$
and close to bound given by
the last entry of Column 2 in Table~1 with $1/\tau(0,0) =10^{7.0458}
\approx 1.1\cdot 10^7$ for $\beta =0$.

\bigskip

Summarizing, the above simple examples confirm the independence of
the theoretical results with respect to the dimension $n$ and the number $m$
of functions in the definition of $f$, as long as the variance
$\|\nabla f(x)-\nabla f_k(x)\|^2$
is bounded independent of $x$ and $n$. The observed results also indicate that the
theoretical results in Table~1 and Table~2 are not overly pessimistic.

\bigskip

\subsubsection{Relaxing the Assumptions}
To test the limits of Algorithm \ref{alg1}
the next example concerns a somewhat more difficult case where again $m=\infty$,
and where each single function $f_k$ carries rather little information
about the function $f$ to be minimized.

\bigskip

For this set of examples a non-singular
randomly generated fixed matrix $A\in\re^{n\times n}$ is chosen
and a scaling factor $\rho :=1/\sqrt n$.
For $k=1,2,\ldots$ random  vectors $r^k$ and $b^k$ are drawn independently
from $n$-variate normal
distributions with expected value $0_n$ and covariance matrices $I_n$ respectively $\rho^2I_n$, i.e. 
$r^k\sim N_n(0_n,I_n)$ and $b^k\sim N_n(0_n,\rho^2I_n)$. The choice of $\rho$ implies that
$E(\|b^k\|_2^2) =1$.
Define $a^k := A r^k$ (so that $a^k\sim N_n(0_n,AA^\T)$) and for given $x\in\re^n$ let
$f_k(x):=\tfrac12((a^k)^\text{T}x)^2+(b^k)^\text{T}x$ with the expected value
\begin{align*}
E\left(f_k(x)\right)
&=  E\left(\frac12((a^k)^\text{T}x)^2+(b^k)^\text{T}x\right)\\
&= \frac12 E\left(x^\text{T}a^k (a^k)^\text{T}x\right)+\underbrace{E\left(b^k\right)^\text{T}}_{=0_n}x\\
&= \frac12x^\text{T} E\left(a^k (a^k)^\text{T}\right)x,
\end{align*}
where the last equation uses the linearity of the expected value for fixed $x$.
Since
\begin{align*}
  E\left( a^k(a^k)^\text{T} \right) = E\left( Ar^k(Ar^k)^\text{T} \right) =
  E\left( Ar^k(r^k)^\text{T}A^\text{T} \right) =  AE\left(r^k(r^k)^\text{T}\right) A^\text{T} = AA^\text{T},
\end{align*}
one can proceed
\begin{align*}
 E\left(f_k(x)\right)
 = \tfrac12x^\text{T}A A^\text{T}x =:f(x).
\end{align*}
And since the fourth momenta of $r^k$ and $b^k$ exist and $f_k$ is a quadratic function of $r^k$ and $b^k$
it follows that $f_k$ has bounded variance and 
$$
\lim_{m\to\infty}\frac1m\sum_{k=1}^mf_k(x)
=f(x)
$$
exists almost surely.
Moreover, $\nabla f(x) = AA^\T x$, which coincides with the expected value of $\nabla f_k(x)$:
\begin{align*}
 E(\nabla f_k(x)) 
 = E\left(a^k(a^k)^\text{T}\right)x+E\left(b^k\right)
= A A^\text{T}x = \nabla f(x).
\end{align*}
The noise defined by $\noise^k:= \nabla f_k(x) - \nabla f(x)=a^k(a^k)^\text{T}x+b^k- A A^\text{T}x$
has expected value
\begin{align*}
 E(\noise^k):= E(\nabla f_k(x)) -\nabla f(x)=  0_n
\end{align*}
and the covariance matrix is given by
\bee \label{variance}
E(\noise^k(\noise^k)^\T)= \rho^2I + AA^\T xx^\T AA^\T + \|A^\T x\|_2^2AA^\T.
\ene
To keep the presentation self-contained a short proof of \mr{variance} is given in the appendix.

If $x\to 0$, the noise terms $\noise^k$ indeed are i.i.d.\ in the limit, as assumed in \mr{iid},
but for larger starting errors $\|x^0\|$ this assumption is violated, more so when
the dimension $n$ grows large.

Note that (since the matrices on the right hand side of \mr{variance} all are positive semidefinite)
$$
\|E(\noise^k(\noise^k)^\T)\|^{1/2} \ge \rho\max \left\{1, \| AA^\T xx^\T AA^\T\|^{1/2}, \sqrt n \,
\|A\,\hbox{Diag}((A^\T x)^2)A^\T\|^{1/2} \right\}
$$
where $\| AA^\T xx^\T AA^\T\|^{1/2}=\|\nabla f(x)\|$. Hence, \mr{variance} implies that
standard deviation of the noise $\noise^k$ (that is added to the gradient
at each iteration of Algorithm \ref{alg1}) always dominates the norm of the gradient itself.

\bigskip

Similar to the analysis in Section \ref{sec.opt}, again, for the numerical experiments
it is sufficient to use a diagonal matrix $D$ in place of $A$. Indeed, consider the
singular value decomposition $A=UDV^\text{T}$ with a
diagonal matrix $D$ containing the singular values of $A$ and
orthogonal matrices $U,V$.
Since $a^k = Ar^k = UDV^\text{T}r^k\sim N_n(0_n,UDDU^\text{T})$
one obtains 
\begin{align*}
 E\left(f_k(x)\right)
=\tfrac12x^\text{T} E\left(a^k (a^k)^\text{T}\right)x
=\tfrac12 x^\text{T}UDDU^\text{T}x
=\tfrac12 z^\text{T}DDz:=\tilde{f}(z)
\end{align*}
with the transformation $z:=U^\text{T}x$.
Likewise, \mr{variance} translates to an equivalent formula for $z$ and $\tilde \noise^k:=U^\T\noise^k$.

The examples in Table~8 and Table~9 refer to
$n=100$, randomly generated $D$
(uniform distribution scaled to $D_{n,n}=n^{-1/2}$ and $D_{1,1}=(10n)^{-1/2}$)
  so that the condition number of $D^2$ is 10.
  The scaling by $n^{-1/2}$ was chosen to compensate for the norm $\|r^k\|$
  which is of the order $n^{1/2}$.
For the choice $k^{max}=1000$ and $\beta=0.7$ (for comparison also $\beta = 0$)
the following average final errors were obtained:

\vbox{
\begin{center}
\begin{tabular}{|c|l||c|c|c|c|c|c|}
  \hline
\multicolumn{2}{|c||}{$\|x^0\|_2$}   & $10^{-1}$ & $10^0$ & $10^1$ & $10^2$ & $10^3$ & $10^4$   \\
  \hline
\multirow{2}{*}{$\|\bar x^{final}\|_2$} & $\beta = 0$  & 0.08 & 0.32 & 3.3 & 33 & 351 & 3218 \\
 &  $\beta = 0.7$ & 0.09 & 0.30 & 2.6 & 24 & 266 & 2817 \\
\hline
\end{tabular}
\vskip 12pt
\centerline{{\bf Table 8:} Final error in dependence of $\|x^0\|$ when $k^{max}=1000$.}
\end{center}
}
When $\|x^0\|_2$ is small, the final iterate is mostly
determined by the noise and it may occur that $\|\bar x^{final}\|_2\gg \|x^0\|_2$.
For the same setting as above and $\|x^0\|_2=10^2$ different values of $k^{max}$ led to
the following results (again $\beta = 0$ is listed for comparison):
\vbox{
\begin{center}
\begin{tabular}{|c|l||c|c|c|c|c|c|}
  \hline
\multicolumn{2}{|c||}{$k^{max}$}   & $10^3$ & $10^4$  & $10^5$ & $10^6$ & $10^7$ & $10^8$ \\
  \hline
\multirow{2}{*}{$\|\bar x^{final}\|_2$} & $\beta = 0$  & 37 & 4.6 & 0.41 & 0.040 & 0.0049 &  6.0e-4   \\
&  $\beta = 0.7$ & 23 & 1.1 & 0.022 & 0.0048 & 0.0018 &  4.6e-4  \\
\hline
\end{tabular}
\vskip 12pt
\centerline{{\bf Table 9:} Final error in dependence of $k^{max}$ when $\|x^0\|_2=10^3$.}
\end{center}
}
Above, it takes about 10 times longer for the algorithm with $\beta = 0$
to reach an error of 0.0048 or 0.0049 than it takes with the choice $\beta = 0.7$.

When multiplying $D$ with a positive constant greater than one, convergence
actually improves (as the noise ratio $\|b^k\|/\|a^k\|$ decreases),
but overall,  
this is an example where the stochastic effects dominate.
This also implies that it may be difficult to improve over the
simple scheme of Algorithm \ref{alg1} for this type of example.

\subsubsection{Non-Quadratic Test Examples}
The intent of this paper of course is to motivate a rule that is more
generally applicable for stochastic gradient descent approaches.
This motivation is supported by the observation that if
the functions $f_i$ are convex and twice continuously differentiable, then
locally, the functions can be closely approximated by quadratic functions for
which the analysis holds.

The MNIST database \cite{deng2012mnist} provides 70000 handwritten digit images with corresponding labels 0,1,\dots,9. Each image consists of $n=28\cdot28=784$ pixels with entries between 0 and 1. To test the performance of Algorithm \ref{alg1}, $m=60000$ of the labeled images are used as training data to find a rule that predicts whether the digits of the remaining 10000 test images are 0 or not. This prediction can be compared with the labels of the test images to determine the false classification rate ``FCR'' i.e.\ the percentage of incorrectly classified images.

For the numerical experiments the data for the $i$-th image is put into the following form:
the vector $a_i\in[0,1]^{784}$ contains the information about the image itself and $b_i\in\{-1,+1\}$
is the corresponding label, where ``+1'' means that the image shows the digit 0 and ``-1'' means
this is not the case.

Let $\sigma:\ \mathbb{R}\to\, ]0,1[\,,\ \sigma(t):= 1/(1+\text{e}^{-t})$ be the logistic function
  with derivative $\sigma'(t)=\sigma(t)(1-\sigma(t))$ and consider $f:\ \mathbb{R}^n\to\mathbb{R}$ defined by
\begin{align*}
 f(x):=-\frac1m\sum_{i=1}^m \log(\sigma(b_i(a_i^\T x))).
\end{align*}
Note that $t\mapsto \log(\sigma(t))$ is a smooth (and asymptotically exact)
approximation of the function $t\mapsto \max\{t,0\}$.
Thus minimizing $f$, approximates the maximization of the terms $\max\{ b_i(a_i^\T x),\ 0\}$ for $1\le i\le m$.
If all terms $b_i(a_i^\T x)$ are positive then 
\bee\label{xstar}
  a_i^\T x>0 \quad\hbox{ for all } i \hbox{ with } b_i=1 \quad\hbox{ and }\quad
  a_i^\T x<0 \quad\hbox{ for all } i \hbox{ with } b_i=-1.
\ene
This motivates the following classification rule:
Let an approximate minimizer $\bar x^{final}$ of $f$ be given and a test image $a_{\text{new}}$. Then $a_{\text{new}}$
is classified to represent the digit ``0'' if $a_{\text{new}}^\T\bar x^{final}>0$, and $a_{\text{new}}$
is classified not to represent the digit ``0'' otherwise.
\\
The derivatives of $f$ are given by
\begin{equation*}
       \begin{aligned}
         \nabla f (x) =  & - \tfrac{1}{m}\sum_{i=1}^{m} b_i a_i \tfrac{ \sigma'(b_i(a_i^\T x))}
                { \sigma(b_i (a_i^\T x))}= - \tfrac{1}{m}\sum_{i=1}^{m} b_i a_i (1-  \sigma(b_i (a_i^\T x))),\\
  \nabla^2 f(x) = & - \tfrac{1}{m}\sum_{i=1}^m b_ia_i (- b_i a_i^\T\sigma'(b_i(a_i^\T x))) \\
  =& \underbrace{\tfrac{1}{m}}_{> 0}\sum_{i=1}^m \underbrace{b_i^2}_{=1} \underbrace{a_ia_i^\T}_{\succeq 0}
  \underbrace{\sigma(b_i(a_i^\T x))(1-\sigma(b_i(a_i^\T x)))}_{> 0} \succeq 0.
 \end{aligned}
\end{equation*}
Here, $f$ is convex since the Hessian matrix is positive semidefinite, i.\,e. $\nabla^2f(x)\succeq 0$.

\bigskip

To estimate a value for the maximum step length $c$ in Algorithm \ref{alg1}, consider the case that
$a_i$ is a random vector whose components are continuously uniformly distributed on $]0,1[\,$. Thus, every component of $a_i$ has expected value $1/2$ and second momentum $1/3$ and the components of $a_i$ are independent of each other. Since
\begin{align*}
  E(a_ia_i^\T)_{r,s} = \left\{
  \begin{tabular}{ll}
   $1/3$ &if $r=s$\\
   $1/4$ &if $r\neq s$
  \end{tabular}
  \right.
\end{align*}
the Hessian matrix can be estimated by
\begin{align*}
 E(\nabla^2 f(0)) = \frac{1}{m}\sum_{i=1}^m E(a_ia_i^\T)\tfrac12(1-\tfrac12)
 = \frac{1}{4m}m \left(\tfrac14\mathds{1}_{n\times n} + \tfrac{1}{12}I_n\right)
 = \tfrac1{16}\mathds{1}_{n\times n} + \tfrac{1}{48}I_n
\end{align*}
with the maximum eigenvalue $n/16+1/48$. This leads to the approximation $c:= 16/n$ used
  for the results in Table~10.

\bigskip

In Table~10 the results of Algorithm \ref{alg1}
with $\alpha=0$ and $\beta\in\{0,0.7\}$
are listed for the initial value $w_0=0_{784}$
and for varying numbers of iterations $k^{max}$:

\vbox{
\begin{center}
\begin{tabular}{|c|l||c|c|c|c|c|}
  \hline
\multicolumn{2}{|c||}{$k^{max}$} & $10^3$ & $10^4$ & $10^5$ & $10^6$ & $10^7$    \\
  \hline
\multirow{2}{*}{$\|\nabla f(\bar x^{final})\|_2$} 
& $\beta=0$ & 6.7e-2 & 1.1e-2 & 2.4e-3 & 1.2e-3 & 7.9e-4\\ \cline{2-7}
&$\beta=0.7$& 5.0e-2 & 8.3e-3 & 1.8e-3 & 7.6e-4 & 4.1e-4\\
\hline
\multirow{2}{*}{FCR} &$\beta=0$& 1.87\% & 0.97\% & 0.74\% & 0.62\% & 0.70\%\\ \cline{2-7}
&$\beta=0.7$ & 1.70\% & 0.96\% & 0.73\% & 0.63\% & 0.70\%\\
 \hline
\end{tabular}
\vskip 12pt
\center{{\bf Table 10:} Norm of the gradient and false classification rate in dependence of $k^{max}$.}
\end{center}
}
For this example it was also possible to apply Newton's method with line search which generated an approximation $x^{opt}$
with $\|\nabla f(x^{opt})\|_2\approx 2.3\cdot 10^{-10}$ and an associated false classification rate of FCR$\approx 0.80\%$.
The Hessian of $f$ at $x^{opt}$ was numerically singular and the spectrum of $\nabla^2 f(x^{opt})$ was
quite dense near zero and did not allow a clear identification of the null space.
It was not possible to give a reliable estimate for the condition number of the Hessian of $f$,
even when restricted to the range space of $\nabla^2 f(x^{opt})$. The Hessian certainly was very far from being
well-conditioned.\\
For all approximate solutions $\bar x^{final}$
generated in Table 10, the distance $\| \bar x^{final}-x^{opt}\|$ was about the same, namely close to 1100.
In particular, a convergence of the iterates of Algorithm~\ref{alg1} to a minimizer of $f$ could not be observed.
Nevertheless a small gain in the rate of convergence with $\beta = 0.7$ compared to $\beta = 0$
could be observed for smaller values of $k^{max}$. For larger values of $k^{max}$
the asymptotic optimality of $\beta = 0$ in \cite{polyak} can be confirmed in the sense that 
there is not much difference of $\beta = 0.7$ and $\beta = 0$.

It is stressed that the intention of this example was not to propose a new classification scheme for MNIST
  -- other classification schemes are certainly better -- but to test Algorithm \ref{alg1} with a somewhat realistic
  example. In particular, the deterioration of the false classification rate in the last column -- and for $x^{opt}$ --
  indicate that the phenomenon of overfitting must be addressed with this approach.

\section{Appendix}

{\bf Proof of \mr{variance}:}\\
Note that
\begin{align*}
E(\noise^k(\noise^k)^\T)&=E((\nabla f_k(x) -\nabla f(x))(\nabla f_k(x) -\nabla f(x))^\T)\\
&=E((aa^\T x+b-AA^\T x)(aa^\T x+b-AA^\T x)^\T)
\end{align*}
where $a = a^k = Ar^k$ and $b=b^k$ with independent normally distributed vectors $r^k\sim N_n(0_n, I)$
and  $b^k\sim N_n(0_n,\rho^2 I)$, a fixed matrix $A$ and a fixed vector $x$.
Multiplying the product in the expected value returns a sum of 9 terms that are considered one by one:\\
The expectation of the constant term $AA^\T xx^\T AA^\T$, of course, is the term itself.\\
The expectations of $-AA^\T xb^\T$ and $-bx^\T AA^\T$ are both zero (since $AA^\T x$ is a fixed vector).\\
The expectation of $bb^\T$ is $\rho^2I$.\\
The expectation of the two terms $-AA^\T xx^\T aa^\T$ and $-aa^\T xx^\T AA^\T$ is $-AA^\T xx^\T AA^\T$ each.\\
The expectation of $aa^\T xx^\T aa^\T$ is given by
$$
E(aa^\T xx^\T aa^\T)=E(A\,rr^\T\underbrace{A^\T xx^\T A}_{=:B}rr^\T A^\T)=A\, E(rr^\T Brr^\T)\, A^\T
$$
where the $i$-th component of $r=r^k\sim N_n(0_n,I_n)$, denoted by $r_i\sim N(0,1)$, has the momenta $E(r_i)=0$, $E(r_i^2)=1$ and $E(r_i^4)=3$. By distinguishing all possible cases one obtains
$$
E(r_ir_jr_pr_q)=\left\{\begin{array}{ll} 3  & \hbox{if} \ i=j=p=q,\\
1  & \hbox{if} \ i=j\not=p=q,\\
1 & \hbox{if} \ i\not =j \ \hbox{and} \ ((i = p,\ j= q) \ \hbox{or} \ (i = q,\ j= p)),\\
0 & \hbox{else}. \end{array}\right.
$$. 
With $i,j\in\{1,\dots,n\}$ the expected value of $rr^\T Brr^\T$ is given componentwise by
\begin{align*}
 E(rr^\T Brr^\T)_{i,j}&=E(e_i^\T rr^\T Brr^\T e_j) = E(r_ir_jr^\T Br) = E\left(r_ir_j\sum_{p,q}B_{p,q}r_pr_q\right)\\
 &=\sum_{p,q}B_{p,q}E(r_ir_jr_pr_q) 
 = \left\{ \begin{matrix}3 B_{i,i}+\sum_{p=1, p\neq i}^nB_{p,p}\ \ \ \text{ if } i=j\\
   2 B_{i,j}\ \ \ \ \ \ \ \ \ \ \ \ \ \ \ \ \ \ \ \ \ \ \ \ \text{ else}\ \ \ \ \ \end{matrix}
 \right.\\
 &=\left\{ \begin{matrix}2 B_{i,i}+\text{tr}(B)\ \ \ \text{ if } i=j\\2 B_{i,j}
   \ \ \ \ \ \ \ \ \ \ \ \ \ \  \text{ else.} \ \ \ \ \end{matrix}
 \right.
\end{align*}
Consequently, one obtains
\begin{align*}
 E(aa^\T xx^\T aa^\T)&=A\, E(rr^\T Brr^\T)\, A^\T\\ &= A\,(2B +\text{tr}(B)I_n )\,A^\T\\
 &= 2AA^\T xx^\T AA^\T + AA^\T\text{tr}(A^\T xx^\T A)\\
 &= 2AA^\T xx^\T AA^\T + \|A^\T x\|_2^2AA^\T.
\end{align*}

Finally, the expectations of
  $bx^\T aa^\T$ and of $aa^\T xb^\T$ both are zero, since $b$ is chosen independently of $a$
and thus also of $x^\T aa^\T$.

\bigskip

\noindent
Summing up  all 9 expectations above leads to \mr{variance}.

\section{Conclusion}
The optimization error and the stochastic error are analyzed for the simple case
  of a stochastic gradient method applied to a
  strongly convex quadratic function  with stochastic gradients that satisfy
  the assumption of being i.i.d. By optimizing both errors for an algorithm with
  finite termination it was possible to modify the known asymptotically optimal parameter selection of
  a stochastic gradient method. As predicted by the analysis in \cite{polyak}, for large numbers
  of iterations $k^{max}$ the results cannot be improved, but for moderate numbers of iterations
  the numerical experiments confirm a gain in accuracy. This gain can be achieved by a simple
  parameter selection and without additional computational cost.
  
  The extension to smooth non-convex functions or to limited memory Quasi-Newton approaches
  such as presented in \cite{burdakov} are the subject of future research.   

\bibliography{literatur}
\bibliographystyle{natdin}

\end{document}